\documentclass{amsart}

\usepackage{graphics,verbatim,color,amsthm}
\usepackage[all]{xy}
\usepackage{hyperref}

\theoremstyle{plain}
\newtheorem{theorem}{Theorem}
\newtheorem{backgroundthm}{Background Theorem}
\newtheorem{proposition}{Proposition}

\newtheorem{lemma}{Lemma}  
 
\newtheorem{definition}{Definition}

\theoremstyle{remark}

\theoremstyle{definition}
\newtheorem{example}{Example}

\def\R{\mathbb{R}}
\def\Z{\mathbb{Z}}

\def\H{\mathbb{H}}

\title{Metric Cones, N-body collisions,   and  Marchal's lemma}

\author{  Richard Montgomery}

\begin{document}

\maketitle

\begin{abstract} Marchal's lemma is the  basic tool for eliminating collisions
when using the direct method of the calculus of  variations to establish   existence of
 ``designer'' solutions to the classical N-body problem.  Our goal here is to understand why Marchal's lemma
holds, by taking a  metric geometry
 perspective  and  employing the   Jacobi-Maupertuis [JM]  metric reformulation of mechanics.  
Using analysis inspired by the  conical metric nature of the  standard Kepler problem at zero-energy, we are able to   manufacture   potentials,
or ``counterexamples'', for which Marchal's lemma
fails.    
These counterexamples,   overlap significantly with results  obtained by  Barutello et al \cite{Terracini}.   A novel feature in our proof for the   counterexample 
is the use of piecewise constant potentials, and the resulting piecewise constant  metrics.  
\end{abstract}

$$ $$

The direct method of the calculus of variations has become  a basic tool for   establishing  the existence of  
interesting  new  periodic solution to    the N-body problem \cite{remarkable}, 
  \cite{Venturelli}.      The possibility that an action minimizer might  suffer a collision  
is the  main theoretical  obstacle to overcome in establishing existence using the direct method.      Marchal's lemma is the fundamental  tool  for overcoming
  this obstacle. 
  
    \begin{backgroundthm} [Marchal's Lemma] [See  \cite{Marchal}, \cite{Chenciner}]  Consider 
      the standard action $\int L dt,  L = K - V$  for any power  law potential $V=V_{\alpha}$ (eq \ref{powerlaw}) of degree $-\alpha$, $\alpha  >0$. 
      (The Newtonian case corresponds to $\alpha = 1$.) 
    Then any   action-minimizer for the fixed endpoint,  fixed-time   problem
   has no interior collision points.
   \end{backgroundthm}
   An `` interior collision point''   means  a collision point  along the solution path which is not one of its
   endpoints.    
   
   Despite this lemma being such a fundamental tool,  I felt like I never  understood
   {\it why} the lemma was true.  In an effort to understand why,  I  will recast the lemma in metric terms,
   and   relate  the validity of the lemma
 to the  {\it inextendibility} of Jacobi-Maupertuis geodesics.    Using this relation, I    construct counterexamples to Marchal's lemma,
 which is to say,  potentials for which the lemma
 fails. These examples give me  a better understanding of why the lemma works.  A central ingredient
 in  constructing these examples  is the conical nature of the JM metric near collision as explained in Proposition \ref{prop_normalform}.
 
 Barutello et al \cite{Terracini}, \cite{Barutello2}    answer  this  ``why?'' question  from  a  somewhat different
 perspective, and in so doing providing quite sharp and illuminating counterexamples to Marchal's lemma.   See  section \ref{BarutelloEtAl} 
 here   for some details on their work and   comparisons to our perspective. 
  
 Both  their work and mine   rely in an essential way on the homogeneity of the potentials.
 It may  be of interest to get rid of this homogeneity condition, perhaps replacing it by a local homogeneity at 
 collision.  In the hopes of doing so we begin without any homogeneity assumptions.

  \section{Set-up}   
  
  Take  configuration space to be  a   Euclidean vector space,  $\R^M$,  endowed with a continuous   function $U$,
  the {\it negative} of the usual potential,
  $$U = -V : \R^M \to (0, \infty ]. $$
   Newton's equations 
  $$\ddot q = \nabla U(q)$$
  are the Euler Lagrange equations for the action whose Lagrangian is
  $$L(q, \dot q)  = K (\dot q) + U(q).$$
  The  conserved energy associated to Newton's  equations  is
  $$H(q, \dot q) = K (\dot q) - U(q)$$
  where 
  $$K(\dot q) = \frac{1}{2} |\dot q|^2$$ is the usual   kinetic energy.  
  
  \begin{definition} By a ``collision point'' we mean a point $p \in \R^M$ for which   $U(p) = \infty$.
  We will also refer to collisions as ``poles''
  \end{definition}
  
  We assume that $U$ is smooth away from the collision points. 
  For simplicity we imagine  that  the collision set is
   a stratified algebraic subvariety although it is not clear how  essential  this is for the   development
   of the theory.

  \begin{example}
[Newtonian N-body  and Power law potentials]  The configuration space of  N point masses moving in $d$-dimensional space
has dimension $M = dN$.   We write  $q_a \in \R^d$ to represent the location of the $a$th body, $a =1, \ldots , N$, so that a vector 
  $q = (q_1, \ldots ,q_N) \in \R^M$    represents the locations of all  N bodies.  
  Write  
    $r_{ab} = | q_a -q_b |$ for  the distance between   body $a$ and body $b$.
  Then the  power law potentials are 
   \begin{equation}
  \label{powerlaw}
   V(q)  = - \kappa \Sigma_{a<b} \frac{m_a m_b}{(r_{ab})^{\alpha}}; \qquad V: (\R^d)^N \to [-\infty, 0).
  \end{equation} 
  Here the $m_a > 0$ represent the masses and $\kappa$ a ``gravitational constant''. 
  The standard Newtonian case corresponds to $\alpha =1$ and $d =3$. 
  A collision occurs in the sense of our definition exactly when a collision occurs
  in the usual sense of   $r_{ab} = 0$ for some $a \ne b$.  
  To get the correct N-body equations we must use  the ``mass inner product''
  $$\langle q , q' \rangle = \Sigma m_a q_a \cdot q'_a$$ to define  the gradient
  in Newton's equations.  Here the dot product of  $q_a \cdot q'_a$ is the standard
  dot product of $\R^d$.  We also must use the mass metric to define kinetic energy.
  \end{example}
   
   By the {\it fixed endpoint, fixed time  action minimization problem} we mean the problem
   where we fix two end points $q_0, q_1 \in \R^M$ and a positive  time $T$ 
   and ask to minimize the action $\int_c L dt$  over all paths $c$ which join
   $q_0$ to $q_1$ in time $T$.
   \begin{definition}  Marchal's lemma holds for the potential $V$
  if every   minimizer for the fixed-end point, fixed time  action minimization problem is free of interior collision points,
   this being true for every choice of endpoints  $q_0, q_1$ and positive time $T$ defining the problem. 
   \end{definition}
   
   \section{The Jacobi-Maupertuis metric}

  Fix a value  $H = h_0$ for the energy $H = K - U$.  Any solution to Newton's
  equations having this  energy will lie  in the associated {\it Hill's region}
   $\{q \in \R^M: h_0 + U(q) \ge 0 \}$. 
   On the Hill region we have the  Jacobi-Maupertuis [JM] metric 
 \begin{equation}
 \label{JMmetric}
 ds^2  = 2(h_0 + U(q)) |dq|^2 
 \end{equation}
 where $|dq|^2$ denotes the standard flat Euclidean metric on $\R^M$. 
  This metric is Riemannian on the open set $0 < h_0 + U < + \infty$ and the well-known
  notion of geodesics and their equations hold here.  In this context the following theorem is well-known.  
  (See for example section 173 of \cite{Wintner}, or exercise 3.4D of \cite{AbrahamMarsden}.) 
   \begin{backgroundthm} Geodesic arcs  for the JM metric at energy $h_0$ which lie in the region
   $0 < h_0 + U < \infty$ are, after reparameterization, solutions to Newton's equations having energy $h_0$.
   Conversely, solution arcs  to Newton's equations which have energy $h_0$ and satisfying
   $0 < h_0 + U < \infty$ are geodesic arcs for the JM metric at energy $h_0$.
   \end{backgroundthm}

   What happens to  geodesics when they pass  through points where the conformal factor $2(h_0 + U)$ vanishes or becomes infinite?
   It is  somewhat unclear what ``geodesic'' even means at such points.
   To give a definition, we start with the Jacobi-Maupertuis length functional $\ell$.  
   Define the length $\ell(c)$  of an absolutely continuous curve 
   $c$ lying in the Hill region to be $\int_c ds$ where   $ds = \sqrt{2(h_0 + U(q))}|dq|$, and  
   allowing $c$ to pass through regions where the conformal factor blows up or vanishes.  
      \begin{definition}  A curve $c$ lying in the Hill region and  joining two points $p$ and $q$
   is a Jacobi-Maupertuis [JM] minimizing geodesic for the energy $h_0$ if $\ell(c) < \infty$
   and if $\ell(c) \le \ell(\sigma)$
   for all other absolutely continuous curves lying in the Hill region and joining $p$ to $q$.
   \end{definition} 
   
    \begin{definition} We say that the potential $V$ satisfies the Jacobi-Maupertuis-Marchal [JM-Marchal] lemma at energy $h_0$,
 if no minimizing JM geodesic for energy $h_0$ has an interior collision point.   \end{definition}

 {\bf Question.} Does the Marchal lemma hold if and only if the JM-Marchal lemma holds for some
 (or every?) energy $h_0$?
 
 We   answer the question in one direction:
 \begin{theorem}
 \label{MarchImpliesJMMarch} If the Marchal lemma holds for a potential $V$, then the zero-energy JM-Marchal lemma holds for this potential $V$.
 \end{theorem}
 
 Under homogeneity assumptions we   obtain the other direction of the  implication.
 \begin{theorem} 
 \label{JMMarchImpliesMarch} Let $V$ be negative,  homogeneous of degree $-\alpha$ 
 for some $\alpha >0$, and smooth away
 from $0$.  If the zero-energy JM-Marchal lemma holds for $V$, 
 then the Marchal lemma holds for $V$.
 \end{theorem}
 
 {\bf Remark.} Under the   assumptions of Theorem  \ref{JMMarchImpliesMarch}   the only collision point is the origin.
 
  {\bf Remark.} The homogeneities  of primary interest in the theorem  are $0 < \alpha < 2$.  This is because the assertion
  of the theorem is   vacuous for $\alpha \ge 2$
  since in this range every path having a collision has infinite action and infinite JM length.

   \begin{theorem} \label{countereg} [Counterexample to Marchal] There are  
  potentials of the form described in Theorem  \ref{JMMarchImpliesMarch} for which the zero energy JM Marchal lemma fails,
  and hence for which, by Theorem \ref{MarchImpliesJMMarch},    Marchal's lemma fails.
  \end{theorem}
  
   See Theorem \ref{JMMarchal} below for a description of the   counterexample potentials.

\section{From Action to JM length.  Proving Theorem \ref{MarchImpliesJMMarch}. }  


For any  real numbers $a, b$ we have
$ab \le \frac{1}{2} (a^2 + b^2)$ with equality if and only if $a-b= 0$.
Setting $a = |\dot q (t)| = \sqrt{2K(\dot q (t))}$ and $b = \sqrt{2U (q(t))}$ we find that
$\sqrt{4K U} \le K + U = L$ with equality if and only if $K-U = 0$.
But the squared norm of the vector $\dot q (t)$ with respect to the zero-energy JM metric $ds^2$  is
$2 U(q(t)) |\dot q (t)|^2 = 4K U$ (see eq (\ref{JMmetric})),  so that  $\sqrt{4K U}dt = ds$
is the integrand for computing the zero-energy JM arclength $\ell(\gamma)$ of a curve $\gamma$.  Integrating our pointwise  inequality yields
$$\ell (\gamma) \le \int_{\gamma} L(\gamma(t), \dot \gamma(t))dt$$
with equality if and only if $H(\gamma(t), \dot \gamma(t))$ is zero a.e.

This inequality holds for any absolutely continuous path $\gamma$ whatsoever.
Apply the inequality  to any curve connecting a point p to a point q in any time interval $[a,b]$
and take the infimum over all such curves.  The left hand side becomes the 0-energy JM distance function.
Since the left-hand side equals the right hand side only
when the curve is parameterized so as to have zero energy, which is to say at speed $|\dot q| = \sqrt{2 U(q)}$, we see that 
the right and left hand side infimums are equal on the set of reparameteriziations of 0-energy JM geodesics
joining p to q.  Note also  that the right hand side is parameterization independent, so the infimum of the
left hand side is also parameterization independent.

We have proved:

\begin{proposition} 
\label{JMiffFTM} A curve $\gamma$  minimizes the zero-energy
JM length among all curves sharing its endpoints if and only if 
$\gamma$, upon being   reparameterized to have zero energy, is a 
{\bf free-time action minimizer}  for the action
among all curves having its same endpoints.
\end{proposition} 

Here we have used
\begin{definition} A free-time minimizer for the fixed endpoint problem for the  action $A(\sigma) = \int_{\sigma} Ldt$ is
a curve $\gamma:[a,b] \to \R^M$ such that for all curves $\sigma:[ c,d] \to \R^M$ for
which $\gamma(a) = \sigma(c),  \gamma(b) = \sigma(d)$ we have
$A(\gamma) \le A(\sigma)$.
\end{definition}
Notice that we allow $d-c \ne b-a$ in the definition.

{\sc Proof of Theorem \ref{MarchImpliesJMMarch}:
Marchal implies JM-Marchal.}

We prove that if the zero-energy Jacobi-Marchal does not hold for the potential $V$,
then the Marchal lemma  fails for $V$.    Suppose, then,  that the zero-energy Jacobi-Marchal
lemma fails for $V$.  Then  
there exist two points $p,q \in \R^M$ and a zero-energy JM-Marchal
minimizer joining them which suffers an interior collision.   Reparameterize this minimizer
so as to have  zero energy \.  By  proposition \ref{JMiffFTM}  this reparameterized  curve is a free-time minimizer
for the action between p and q. Free-time minimizers are automatically   fixed time minimizers, the time $T$
being the total time needed to connect the two points in the reparameterization \footnote{The question arises: is $T$ finite? The answer is yes.
To see this we can suppose that the JM minimizer is parameterized by arclength $s$, so that $\sqrt{2 U(q(s))}|dq/ds | =1$.
The new parameterization variable  $\tau$ is to satisfy $\frac{1}{2} |dq / d \tau |^2 - U(q(\tau)) = 0$.  The relation  $dq /d \tau = (dq/ds) (ds/ d \tau)$
and some algebra yields  $d \tau = ds / 2 U (q(s))$.  Since $U \to \infty$ at collision,  the latter factor is integrable with finite integral.
Indeed we could have begun by taking the initial endpoints p and q sufficiently close to collision so that all along the JM minimizer 
joining them we have $U \ge 1$. }.  Since the   minimizer has
an interior   collision point,   Marchal's lemma fails for $V$.  QED.

\section{The Kepler Cone}
\label{KeplerCone}

Our metric understanding of Marchal's lemma began by working out  details in the case of the Kepler problem
in the plane.  The    potential   is $V (q) = - 1/ |q|$, so 
   $U = + 1/|q|$ and  the Jacobi-Maupertuis metric at energy 0   is 
  $$ds^2  = \frac{2}{|q|} |dq|^2$$
 where $|dq|^2 = dx^2 + dy^2$ is the standard Euclidean metric. 
 In polar coordinates $r, \theta$ we have  $|q| = r,  |dq|^2 = dr^2 + r^2 d \theta ^2$, 
 so 
  $$ds^2 _{JM} = \frac{2}{r} (dr^2 + r^2 d \theta ^2) = 2 (\frac{d r^2}{r} + r d \theta ^2).$$
Two metrics related to each other by a positive
 constant have the same geodesics,  so we can delete the overall   factor $2$
 and work  with the metric
 $$ds^2 = \frac{d r^2}{r} + r d \theta ^2$$
 on the plane.

 We  will  now  put this metric  into ``standard  conical normal form'': 
 $d \rho^2 + c ^2  \rho^2  d \theta^2$,  where $c$ is a constant.  We do so  by  a change of variables
 $\rho = \rho (r), \theta = \theta$.  To find $\rho$ set  
 $d \rho^2 = \frac{d r^2}{r}$ or $d \rho = \frac{d r} {r^{1/2}}$  and  integrate, using $\rho(0) = 0$.
  We obtain
 \begin{equation}
 \label{changevars1}
 \rho = 2 r^{1/2} \text{ with inverse }  r = (\frac{\rho}{2})^2, 
 \end{equation}
 and the desired form
 $$ ds^2 = d \rho ^2 + (\frac{1}{2})^2 \rho^2 d \theta^2,$$
 which is  the   metric   of   {\it the  cone over a circle of radius} $c= 1/2$.
 
 {\sc Cones over Circles.}
 
 Let $\theta$ be an angular coordinate, $c$ a positive constant,  and form the metric
 $$ds^2 = d\rho^2 + c^2 \rho d \theta^2, \rho \ge 0, 0 \le \theta \le 2 \pi.$$
 The change of variables
  $$\phi = c \theta$$
  changes this metric into  
  $d \rho ^2 + \rho^2 d \phi^2$
  which is the  Euclidean metric on the plane, written in polar coordinates $(\rho, \phi)$, except now  $\phi$
  is subject to the constraints  
  $0 \le \phi \le 2 \pi c$. For  $c < 1$ these constraints   define a  sector in the plane
    bounded by the rays   $\phi = 0$ and  $\phi = 2 \pi c$.  
  Since $\theta = 0$ and $\theta = 2 \pi$ are identified, we must 
  glue  the two  bounding rays of this sector together   to   form our cone.  The result is a standard cone made by gluing.
  We call the  $(\rho, \phi) $ sectorial representation of the cone the ``cut flattened cone''.
  
  Our particular case $c =1/2$  for Kepler corresponds to a half-plane, with  the bounding rays
  part of a single line, the line which  forms the half-plane's boundary.  
  This half space is the fundamental domain for the action of $\Z_2$ on the Euclidean plane
  by $(x,y) \mapsto (-x, -y)$, which, on the boundary of the half-plane corresponds to the identification of the two
  bounding rays leaving the origin.  Consequently, the Kepler  cone  is isometric to  the metric quotient
  $\R^2/ \Z_2$.
  
   {\bf  Remark.}  It has been known for some time that   the Jacobi-Maupertuis metric for  the zero energy Kepler problem  is flat
   away from collisions.   See
   for example  section 244  of \cite{Wintner}, where  Winter   remarks that the Gaussian curvature of the  JM metric at energy $h_0$
 is $-\frac{h_0}{4} [r (h_0 + 1/r)]^{-3}$   in the open  region $0 < (h_0 + 1/r) < + \infty$.   
  
 Minimizing geodesics  on the cut flattened cone are Euclidean line segments.
      Using  this fact  we  can prove: 
 \begin{lemma} On the cone over a circle of radius $c$, $c < 1$, any  geodesic ending at the  cone point is inextendible:
 it cannot be extended and remain geodesic. 
 \end{lemma}  
 We refer the reader to our Appendix on Metric Geometry for the precise notion of `geodesic' and `inextendible geodesic'
 in a general metric space.  See in particular definition \ref{inextendible} of that Appendix. Our definition of the JM-Marchal lemma holding
 is equivalent to the assertion that the every geodesic ending in collision is inextendible. 
 
  \begin{figure}[h]
\scalebox{0.4}{\includegraphics{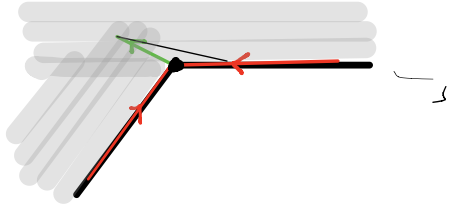}}
\caption{The conical metric can be flattened to a sector, whose bounding rays are glued to form the cone.
An incoming geodesic can be rotated to form this glued boundary.  
The depicted  shortening   shows that no matter how we geodesically extend
the incoming ray, we can shorten the result while avoiding collision with the  cone point . }
\label{fig_sector}
\end{figure}
  {\bf Proof.}   Geodesics are line segments.
  We must show that the concatenation of a line segment coming in to the cone point with one
  exiting the cone point always  
  fails to minimize length. So consider an incoming geodesic ray (in red) and outgoing geodesic segment (green).  
   See figure \ref{fig_sector}.  By rotational symmetry, we can assume that  the 
  incoming ray   coincides with
 one of the sector's bounding rays.    Because the two bounding rays are identified to form the cone, the incoming ray  actually    
  corresponds to {\bf both} bounding rays.  Since $c < 1$, however we orient the   outgoing geodesic segment (green), its angle with one or the other of the two
 bounding rays is less than $\pi$. Hence we can ``cut the corner''  (in black in the figure),  skipping the cone point,  and shortening  the
 resulting concatenated curve.  
 QED

Lemma 1 establishes the validity of the zero-energy Jacobi-Marchal lemma for the Kepler case, and hence, upon invoking Theorem \ref{JMMarchImpliesMarch}, 
 the Marchal lemma. 

 {\sc  Metrics  Cones, Generally.} 
 
  A brief  discussion of more general  metrics cones is in order.
    

If $Y$ is a manifold with Riemannian metric $ds^2 _Y$, we  form the
``metric cone over $Y$ ''  by putting the
metric 
$$d \rho ^2 + \rho ^2 ds^2 _Y ,  \rho > 0  $$
on $(0, \infty) \times Y$, and noting that as $\rho \to 0$
the metric on the $Y$ factor shrinks to zero.
Thus, we crunch $0 \times Y$ to a single point,  called the cone point.
Topologically,  crunching is achieved by dividing out $[0, \infty) \times Y$
by the equivalence relation  $\sim$  in which all points $(0, y)$ are identified with each other.
The resulting topological space  $([0, \infty) \times Y) / \sim$
is   the ``cone over $Y$'',  denoted  $Cone(Y)$.  The function $\rho$ is then the distance
from the cone point, and the `spheres'' $\rho = \rho_0$  about the cone point are copies of $Y$, with the
$Y$'s metric
scaled by $\rho_0^2$.     The   case of a cone over the circle corresponds to the case
where $Y$ is a circle of radius $c$.
 For the theory of general length spaces, and for cones over them, we refer the reader to Burago et al \cite{Burago}.

   \section{Cones and the  zero energy JM metric} 
 
 The zero energy JM   metric is 
$$ds^2 = 2 U(q) |dq|^2$$
where $|dq|^2$ is the standard Euclidean metric on $\R^M$.  
Use spherical coordinates $(r, s) \in [0, \infty) \times S^{n-1}$
with $q = r s \in \R^n$ , $r  = |q|$,   $s$ a unit vector. 
{\it Now assume that $U$ is homogeneous of degree $-\alpha$.}
Then we have the   ``shape  potential''   $\hat U$, or ``normalized potential''  
\begin{equation} \hat U: S^{n-1} \to (0, \infty] \text{ defined by   }
U (r s) = r^{-\alpha} \hat U (s),  \| s \| = 1.
\label{normalized}
\end{equation}
The kinetic energy in spherical coordinates is  
$$|dq|^2 = dr ^2 +  r^2 ds ^2_{sphere}  $$
where $ds ^2_{sphere}$ is the standard round metric on
the  unit sphere $S^{n-1}$.   Thus
$$U ds^2 _K = r^{-\alpha} \hat U (dr^2 + r^2 ds^2 _{sphere})  = \hat U (s) (r ^{ -\alpha}dr^2 + r^{2 - \alpha} ds^2 _{sphere}).$$
Solve 
$$ d \rho ^2 = r^{-\alpha} dr ^2$$
with the boundary condition  $\rho =0$ when $r = 0$ to get, for $0 < \alpha < 2$ \footnote{If $\alpha \ge 2$ then the integral of $r^{-\alpha/2} dr$ from $\epsilon$ to $1$ diverges as $\epsilon \to 0$, this
change of variables cannot be made.}
$$\rho =   \frac{1}{c( \alpha)} r^{1 - \alpha/2},   \qquad r^{2 -\alpha}  =  c(\alpha)^2 \rho ^2,  \qquad c(\alpha)  = \frac{2 - \alpha}{2} $$
and finally 
\begin{equation}
\label{NormalFormA}
ds ^2  =  \hat U (s)  ( d \rho ^2 +  c^2 \rho^2 ds^2 _{sphere} ),  c(\alpha)  = \frac{2 - \alpha}{2}. 
\end{equation}
Observe that $d \rho ^2 +  c^2 \rho^2 ds^2 _{sphere}$ is the metric for the cone over the sphere of radius $c = c(\alpha)$. 
We have shown
\begin{proposition}
\label{prop_normalform} The zero-energy JM metric for any negative potential which is homogeneous
of degree $-\alpha$,  $0 < \alpha < 2$ is  given by the expression  (\ref{NormalFormA}),
which is that of a metric conformal to 
 the cone over the sphere of radius $c( \alpha) <1$ with  conformal factor
 the normalized shape potential $\hat U$.  
\end{proposition}

\section {The counterexample: Theorem \ref{countereg}.}

\begin{figure}[h]
\scalebox{0.4}
{\includegraphics{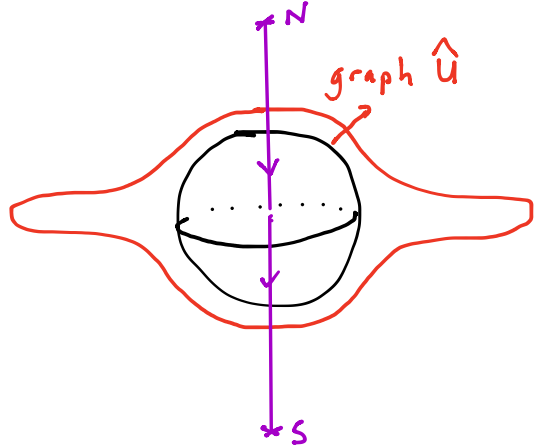}}
\caption{The shape of the normalized potential for the counterexamples.
The minimizer buries right through the center of the sphere, passing through collision
as it travels between the North and South poles.  }
\label{fig_potential}
\end{figure}

The normal form,  eq (\ref{NormalFormA}) described in Proposition \ref{prop_normalform}, provides the idea of   how to construct a counterexample
to Marchal's lemma, i.e. a potential for which the lemma fails.
Design  $\hat U$ to have  absolute minima at   the poles N, S of the sphere, while at the same time  to be  very large in a   band surrounding    the equator.
Look for minimizers from $p=N$ to $q=S$.  
Burying straight down through the earth by travelling the Euclidean line segment in $\R^M$  joining   N to S is a collision path which
 will be   much shorter than travelling along the earth's surface   $\rho =1$  from N to S since in so doing we
 must climb the high mountains surrounding the equatorial band.  
See figure \ref{fig_potential}.
With  work this simple idea can be promoted to a proof of Theorem \ref{countereg}.

Let  
$$z(s)   = s_M:=  \langle s , e_M \rangle$$
 denote the   height coordinate  of a point $s$  on the sphere so that the North and South poles of the sphere, N and S,  are
given by $z = \pm 1$ while its  equator (an $M-2$-sphere) is   defined by $z = 0$.  (Here $e_M$ is the last basis
vector of our Euclidean configuration space $\R^M$.)
For $0 < \delta < 1$ the locus   $|z| < \delta$ is an equatorial band of thickness $2 sin^{-1}(\delta/2)$.  
Suppose   there are  positive constants $\epsilon < 1$,  and   $m < M$ such that  $\hat U$ satisfies  
\begin{itemize}
\label{conditions}
\item{(A)} $\hat U$ achieves its absolute minimum value of $m$ at the two poles N and S 
\item{(B)}  $\hat U   \ge M$ for $|z|   \le \delta$. 
\end{itemize}   
\begin{theorem}
\label{JMMarchal} Suppose that the normalized potential satisfies (A) and (B) above, that 
  the degree of homogeneity of the potential is   $-\alpha$ with  $0 < \alpha < 2$, 
and   that $M \delta  \ge \frac{m}{c(\alpha)}$, where $c(\alpha) = \frac{2-\alpha}{2}$.
Identify the sphere with the locus $\rho = 1$ as per
the representation of eq. (\ref{NormalFormA}) and proposition \ref{prop_normalform}.  
Then the   JM-minimizing geodesic connecting the two poles N and S   is the Euclidean line segment.
In  particular, this  minimizer passes through collision,showing that the JM Marchal lemma fails for this potential.  
\end{theorem}

Theorem \ref{countereg} is simply a restatement of Theorem \ref{JMMarchal},   combined with Theorem \ref{MarchImpliesJMMarch}.

\subsection{Proof of Theorem \ref{JMMarchal}.}   It follows   immediately  from condition  (A)   that the minimizing geodesics
connecting  the sphere $\rho =1$ to total collision at the  origin  $\rho =0$  are the   Euclidean line segments obtained by fixing the shape 
$s$ to be one of the two poles and  letting $\rho$ vary from $1$ to $0$. 
The length of either segment is $m$, so their concatenation, the line segment  $\gamma_*$ 
described in the statement of the theorem,   has length   $2m$, 
connects $N$ to $S$,  and has  $0$ as an interior collision point. 
We must show 
 that any   curve $\gamma$ joining $N$ to $S$ and avoiding collision has length
greater than $2m$.  

Replace $\hat U$ by the piecewise constant function
\begin{equation}\hat U_{step} (s) = \begin{cases} m,  |z(s)| >  \delta , \\
M ,  |z(s)| \le \delta
\end{cases}
\label{step1}
\end{equation}
with corresponding piecewise conical
 metric 
 \begin{equation}
 ds^2 _{step} =  \hat U_{step} (s)  ( d \rho ^2 +  c^2 \rho ^2   ds^2 _{sphere}   ) 
 \label{step2}
 \end{equation} 
Since $\hat U_{step} \le \hat U$ and since the two functions agree at the poles
we have $\ell_{step} (\gamma) \le \ell (\gamma)$ while  $\ell_{step} (\gamma_*) = \ell (\gamma_*)$.   Thus, to prove the theorem  it suffices to show that
  $\ell_{step}(\gamma) >  2 m$ for  $\gamma$   any curve joining $N$ to $S$ and not passing through the origin.

Suppose  $\gamma$ is such a curve.  Since $\gamma$ joins N to S and    $\rho>0$ along $\gamma$,
the spherical part $s(t)$ of $\gamma(t)  = (\rho(t), s(t))$ must  cross the equator at some point
E.  Consider the Euclidean half plane  $\H \subset \R^M$  spanned by N and E, and lying on the E side of the line NS.
This half-plane is parameterized by $\rho$ and an angle $\phi$,  with $0 \le \phi \le \pi$ which parameterizes the semi-circular
 longitude  N E S.
We first show  that  we may assume that $\gamma$ lies
inside this half plane.  Consider the projection operator 
$$pr: \R^M \to \H$$
onto the half plane which can be obtained by   rotating  the 
spherical part  $s$  of the  point with coordinates $(\rho, s)$  until it lies on the longitude NES, while  keeping the radial part $\rho$ fixed.
Write    $\omega_{\perp} \in S^{M-2} $ for a point lying on the equator which is   the  unit $M-2$-dimensional Euclidean 
sphere  orthogonal to $N$.  Then  any point $s$ of the sphere can be written
$$s = \cos(\phi ) N + \sin (\phi) \omega_{\perp} \text{ so that  } pr(s) = \cos(\phi) N + \sin(\phi) E.$$
Thus    $pr( \rho,  s) = (\rho,  pr(s))$.
  To see that $\ell_{step}(pr(\gamma)) \le \ell_{step}(\gamma)$
  observe that we can write the spherical element of arclength occuring in eq (\ref{NormalFormA})
  as $ds^2 _{sphere} = d \phi ^2 + \sin^2 \phi d \omega_{\perp} ^ 2$ 
  and $d \omega_{\perp}^2 = 0$ along $pr(\gamma)$ while
  $\hat U_{step}(pr(s)) = \hat U_{step}(s)$.  It follows that
 along the curve $pr(\gamma)$ 
  the integrand $ds_{step}$ used to compute   $\ell_{step}$ is less than  the same  integrand  along for $\gamma$
  and 
  consequently $\ell_{step}(pr(\gamma)) \le \ell_{step}(\gamma)$ with equality if and only
  if $\gamma = pr(\gamma)$.

\begin{figure}[h]
\scalebox{0.4}{\includegraphics{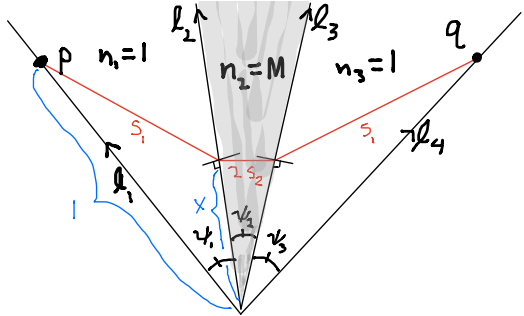}}
\caption{The metric becomes piecewise Euclidean, so geodesics (red) are p.w. linear curves subject to Snell's law.
The $n_i$ indicate indices of refraction, which is to say, constant scale factors with which to multiply the standard
Euclidean metric. }
\label{fig_Snell}
\end{figure}

  We have reduced our problem to a problem of computing lengths  for curves on the half plane
  relative to the step metric. 
  The  metric on the half plane has the
form
 \begin{equation}
 ds^2 _{step} =  \hat U_{step} (\phi)  ( d \rho ^2 +  c^2 \rho ^2   d \phi^2 ) 
 \label{step2}
 \end{equation} 
where we have used polar coordinates   $\rho, \phi$, $0 \le \rho \le \infty, 0 \le \phi \le \pi$ to coordinatize the half-plane. 
This half plane is cut into three sectors, namely $0 \le \phi \le \pi/2 - \delta,
\pi/2 - \delta < \phi < \pi/2 + \delta$ and $\pi/2 + \delta < \phi < \pi$ on which $\hat U_{step}$ 
is constant.  Now any metric of the form  $A( d r ^2 +  c^2 r ^2   d \phi^2 ) $, $A$ a positive constant, 
is locally  isometric to the flat  Euclidean metric, as we saw earlier in the Kepler section by using the  trick of making the substitution  $\theta = c \phi$
which converts the metric to $A (d\rho ^2 + \rho^2 d \theta^2)$. 
At this stage it   becomes crucial   that  $c < 1$, as this factor  shrinks the opening angle  $\pi$ of the half plane to the  
 angle of $c \pi < \pi$.     Writing 
$$\delta = sin(\phi_* /2),$$ the  opening angles of our three sectors
are $\psi_1 = c (\pi - \phi_* )/2,  \psi_2 =  c \phi_*$, and $\psi_3 = c (\pi - \phi_*)/2$.  See figure \ref{fig_Snell}

Our step  metric on the half-plane  consists of a Euclidean metric on each   sector, with  jump discontinuity as we pass from one sector  to the other. Geodesics
within each sector are  Euclidean lines.  The direction of the line
suffers a jump discontinuity in crossing from one sector to another according to  Snell's law of optics. The constants $m$ and $M$  
play the role of the index of refraction in Snell's law.

We  now give the final piece of  intuition underlying our proof.   
If  $M =m$ then the metric is a uniform  Euclidean metric across all 
sectors and so  geodesics are Euclidean straight lines, rather than piecewise linear curves.  
The metric on the half plane has total   opening angle  $c \pi < \pi$ between its bounding rays,  
rays $0N$ and $0S$, so that  the line segment  joining    N to  S    does not   intersect  $0$.
Now fix  $m$ and let  $M$  increase.  The   minimizing geodesic  bends closer to the vertex,
 spending less time in the middle region.   See the red curve in figure \ref{fig_Snell}. 
This process is monotonic in M. Eventually,  at some critical value $M= M_c$ the minimizer from N to S disappears into the
bounding rays. From then on, for all $M > M_c$,  the minimizer coincides with $\gamma_*$.
The proof is finished by showing that $M_c < \frac{m}{c(\alpha) \delta}$.

We now fill in the details.  
Label  the bounding rays of the three sectors  
 $\ell_1, \ell_2, \ell_3, \ell_4$
 so that N lies on  
$\ell_1$, $S$ lies on $\ell_4$ and  $\ell_2$ and $\ell_3$ bound the central sector whose
index of  refraction is M.  
Reflection  $z \mapsto -z$ leaves $\hat U_{step}$ invariant and so   is an
isometry of the  step metric.  Since this reflection takes N to S,
it follows that any minimizer from N to S is invariant under this reflection.
In particular, if the minimizer does enter into the interior of the middle sector, crossing ray $\ell_2$
a distance $x$ from the vertex, then it must leave that  middle sector along $\ell_3$ at  the same 
distance $x$  from the vertex.     Here we  measure   $x$  
relative to the   underlying Euclidean metric
$d\rho ^2 + \rho^2 d \theta^2$, which is the metric we  multiply by $M$ or $m$ to get the metric in the various sectors.
  We have reduced the proof of the proposition   to a single variable calculus
  problem, that of minimizing the step  lengths  
  of our  one parameter family of  ``test curves'' .  Again, see figure  \ref{fig_Snell}. 
  
Let us write $S(x)$  
for the step-length of the test curve  labelled by $x$.
Then $S(x) =  m s_1 (x) + M s_2 (x)  + ms_3 (x) = 2m s_1 (x)  + Ms_2 (x)$
where the $s_i$ are the lengths of the line segments of the test curve   in the underlying  Euclidean metric
and we use the reflectional symmetry to get that $s_1 (x) = s_3 (x)$. 
By the   law of cosines  
$$s_1^2 = 1  + x ^2 -2x \cos(\psi_1),$$
while, from trigonometry   
$$s_2 =  x sin(\psi_2/2).$$ 
Thus 
$$S (x) =  m \sqrt{ 1 + x^2 - 2x \cos(\psi_1)} + M x \sin(\psi_2/2).$$
Differentiating  with respect to $x$ yields
$$ \frac{d S}{dx}  = m \frac{x -  \cos (\psi_1)}{s_1}  + M \sin(\psi_2 /2)$$
We will show that this derivative is always positive for $x > 0$, provided
the condition $M \ge m / (c \delta)$ holds.  

Set $h = x -  \cos (\psi_1)$ and observe that $s_1 = \sqrt{ 1 + h ^2 - \cos^2 (\psi_1)} = \sqrt{ \sin^2 (\psi_1) + h^2}$.
Consequently, the first term of the derivative $d  S/dx$ is  $m \frac{x -  \cos (\psi_1)}{s_1}  =m  \frac{h}{\sqrt{h^2 + \sin^2 (\psi_1)}}$.
which  is always less than or equal to $m$ in absolute value.  Thus  $M \sin(\psi_2/2) > m$
implies  that $d S/dx > 0$.  Now $\psi_2/ 2 = c \phi_* /2$ and $\sin(\phi_* /2) = \delta$ where
we recall the significance of $\delta$ was that $\hat U > M$ for $|z| \le \delta$.   See item (B)
of the conditions \ref{conditions} above.  Now use the inequality $\sin(c \theta)  >  c \sin(\theta)$  which is valid for $0 < |c \theta| < \pi$ and $0 < c < 1$.
We see that $M \sin(\psi_2/ 2) > c M \delta$.  Thus $cM \delta \ge m$ implies that
$d S /dx > 0$ for $x > 0$.    The positivity of this derivative  means that the
length of these test curves increase monotonically from their absolute  minimum value of $2m$ when $x = 0$,
and thus there is  no minimizer  interior to the sector.  QED.

\section{}

 \section{Homogeneity. Blow-up.  Reduction to zero energy. Theorem \ref{JMMarchImpliesMarch}}
 \label{sec:blowup}
 
  Throughout this section we assume   that $V$ is negative,  homogeneous
  of degree $-\alpha$, $0 < \alpha < 2$ and smooth away from $0$. 
  We will begin using blow-up to show how the   zero energy and nonzero energy JM Marchal 
  lemma are related. 
  
  A key step in the usual  proof of the Marchal lemma is blow-up: a
  rescaling argument which reduces the investigation of  action minimizers with
  collision to the case of  zero energy action minimizers having a collision.  Blow-up is based on the fact that 
   if $U$ is homogeneous of degree $-\alpha$ and 
  if  $q(t)$ solves  the corresponding Newton's equations   then $q_{\lambda} (t) = \lambda q (\lambda ^{-\nu(\alpha)} t)$
with 
  $\nu = (1 + (\alpha)/2))$ also solves Newton's equations.  If $q(t)$ had energy $H$ then $q_{\lambda} (t)$ has energy $\lambda^{-\alpha} H$.  
  
  We proceed to   a metric version of blow-up.  Consider the dilation map $F_{\lambda} (q)  =  \lambda q$
  of $\R^M$.  Pull back  the JM metric $ds^2 _H= 2(H+ U(q)) |dq|^2$ by $F_{\lambda}$ to obtain
  $F_{\lambda} ^* ds^2 _H = 2(H + \lambda^{-\alpha} U(q) \lambda ^2 |dq|^2 = \lambda ^{2-\alpha} (\lambda^{\alpha} H + U(q)) |dq|^2$.
  Thus 
  $$F_{\lambda}^* ds^2 _H = \lambda ^{2-\alpha} ds^2 _{\lambda^{\alpha} H}.$$
  Now a metric  and a constant times that metric have the same geodesics, and if
  a geodesic is a minimizer for one, then it is a minimizer for the other.   It follows that if $c$ is
  a    minimizing  geodesic for $ds^2_H$ joining $A$ to $0$
  then  $F_{\lambda} ^{-1}(c)$ is a minimizing geodesic for $ds^2 _{\lambda^{\alpha} H}$
  joining  $\frac{1}{\lambda} A$ to   $0$.  As $\lambda \to 0$ the  curves $F_{\lambda} ^{-1}(c)$
  are uniformly bounded on compact sets , so a subsequence of them converges to a
   minimizing  geodesic for the zero-energy JM metric  $ds^2 _0$, one  which  ends at the cone point.
   Similarly, if for some $H \ne 0$ there is a minimizing geodesic which ends at the cone point
   and can be extended past it, then by dilating and taking subsequences, we arrive at extendible minimizing
   geodesics passing through the cone point.
   
   We have proved
   \begin{proposition}  Suppose that $V$ is homogeneous of degree $-\alpha$
   for $0 < \alpha < 2$, negative, and smooth away from zero. Then the  JM-Marchal lemma holds for   $V$ at energy   $h_0$
   if and only if it holds for $V$ at energy $0$.
   \end{proposition}

\subsection{Proof of theorem \ref{JMMarchImpliesMarch} }
[due to  Andrea Venturelli]
We proceed by proving the contrapositive.  
Suppose that the Marchal lemma for total collisions \footnote{i.e. the only collision available under our hypothesis, $q=0$}  fails for a particular potential of homogeneity
$\alpha$.   We will show that the JM-Marchall lemma also fails for this potential.   
The  failure of the Marchal lemma for total collisions means that    there 
exists  a curve  $\gamma(t)$  which has an internal total collision and is a fixed time minimizer between its
endpoints.  Translate time so   the total collision occurs at time $t=0$.  
By a standard argument (\cite{Venturelli}, \cite{Chenciner}, section 3.2.1),  the rescaled  family 
 $\lambda^{\nu} \gamma(\lambda t)$, with  $\nu = (1 + (\alpha)/2))$,   converges  as $\lambda \to 0$,  to 
a curve  $$\gamma_* = \gamma_+ * \gamma_-$$ which is the 
concatenation of two parabolic homothetic solutions,   $\gamma_- $ defined for 
$t \le 0$ and $\gamma_+$   for $t \ge 0$.  (The convergence is uniform on bounded intervals containing $0$.)
Moreover,  $\gamma_*$ is a global {\bf fixed time} minimzer: that is,  for each   pair  of times $a<  b $ the
segment 
$\gamma_* ([a, b])$ is a {\bf fixed time}   action minimizers between its  endpoints
$\gamma_* (a)$ and $\gamma_* (b)$.  If we knew the concatenation
was a {\it free time minimizer}, rather than just a fixed-time minimizer,
 between all pairs of its points then we would be done, by Proposition \ref{JMiffFTM} .
 But we don't know that yet.  

We now argue by contradiction.  Suppose that the concatenation $\gamma_*$  is not a  global free time minimizer.
Then there must exist two points  along $\gamma_*$ , one before collision, one after,
for which there is curve $\eta$  joining these two points and having smaller action.
Write the points as  $\gamma_- (\alpha)$  and $\gamma_+ (\beta)$,
for  $\alpha  < 0 < \beta$.  Then the   curve  $\eta: [a, b]  \to \R^M$  satisfies  
$\gamma_- (\alpha) = \eta (a),  \gamma_+ (\beta) = \eta (b)$
  $A(\eta) < A(\gamma_* |_{[\alpha, \beta]})$ and $b-a \ne  \beta-\alpha$.
  Set $c = (\beta - \alpha) - (b-a)$.  For large  positive $T$ consider the concatentation
   $$ y = \gamma_+ |_ {[T, \beta] } * \eta * \gamma_- |_{ [\alpha, - T] } ,$$  with   two  of the curves in this concatenation requiring
  a time shift in their parameterizations to guarantee that they take off when the previous curve ends.
  The action of $A(y)$ is
  $$A(y) = A(\gamma_* |_{[-T,T]}) -  A(\gamma_* |_{[\alpha, \beta]}) +  A(\eta) = A_0 - \delta > A_0: = A(\gamma_* |_{[-T,T]})$$ 
  where 
  $$\delta = A(\gamma_{[\alpha, \beta]}) - A (\eta) > 0.$$
  However, the  the interval parameterizing   $y$ is  not  $[-T, T]$
  but rather $[-T, T+ c]$ whose length is   $2T+c$.  
  To finish off the proof we  must reparameterize $y$ so as to be parameterized an interval of length $2T$  
  with the penalty of possibly increasing  the  action, but not enough to swamp the $- \delta$.
  We will estimate that this  this increase due to reparameterization is   $O(1/T)$, so that  taking $T$ sufficiently large will 
  complete the proof.
  
  The ratio of the two intervals of parameterization is  
  $$ \lambda = \frac{2T+ c}{2T} = 1 + \frac{c}{2T}.$$
  So set  
  $$x(t) = y(\lambda t)$$
  yielding a  reparameterization $x$  of $y$ by an interval of length $2T$.  
If the  action $A (y) = \int (K + U )dt$ then
 one computes that   $A(x) =  \lambda \int  K dt  + \frac{1}{\lambda}  \int U dt \le M A(y) $
  where  $M = max{\lambda,  \lambda ^{-1} } $. 
  But $A(y) = A_0 - \delta$ where $A_0 = A(\gamma_*  |_{[-T,T]} )$,
  and $M  = 1 + O(1/T)$.  Thus
  $A(x) = (1 + O(1/T)) (A_0 - \delta ) = A_0 - \delta + O(1/T)$
  which is less than $A_0$ and gives us our contradiction.
  
  QED.

\section{Comparison with the Work of Barutello et al}
\label{BarutelloEtAl}

In a tour-de-force of variational and dynamical  analysis Barutello, Terracini and Verzini  \cite{Terracini}, \cite{Barutello2}   have thoroughly investigated
a class of problems very  similar to ours.   They consider negatively homogeneous
 potentials $V(r s) = - r^{-\alpha} \hat U (s)$ for which $\hat U$ is positive,  sufficiently smooth and takes
 on its minimum   at precisely two distinct points $\xi_-, \xi_+$.   
  The   authors also
assume    both  minima are  nondegenerate.  
 Denote the set of all such normalized potentials on the sphere by   ${\mathcal P} = {\mathcal P} (\xi_-, \xi_+) \subset C^2 (S^{M-1})$.  Then they coordinatize
 the space of all their potentials $V$  by ${\mathcal P} \times (0, 2)$ and where $\alpha \in (0,2)$.
 For each such potential they ask ``{\bf does the free time action minimizer joining $\xi_-$ to $\xi_+$ pass through
 total collision ($q=0$)}?''  If the answer is `yes' they call the potential labelled by $(\hat U, \alpha)$ an `IN' potential, and otherwise they call
 that potential  an ``OUT''
 potential, in this manner decomposing the space of all potentials into two disjoint sets.    Their main result is that there is a continuous function  $f: {\mathcal P} \to (0,2)$ such that
 if $\alpha \le f(\hat U)$ then  $(\hat U, \alpha)$  is an IN potential, while if $\alpha > f(\hat U)$ then this  potential
 is an OUT potential.   (The function $f(\hat U)$ is denoted $\gamma(\hat U, 0^+)$ in their paper.)   
 
 From the perspective of Barutello et al  then, our  theorem \ref{JMMarchal}  asserts that if  
   $M \delta > m/c(\alpha)$ then the  potential is of  ``IN'' type. 
 Doing some algebra, we see that this inequality holds if and only if $2( 1 - \frac{m}{M \delta}) > \alpha$,
 which logically implies the   estimate    $f(\hat U)  \ge 2( 1 -  \frac{m}{M \delta})$.  
  
Barutello et al    have a quite pleasing characterization of the  value $\alpha_*= f(\hat U)$
as a ``phase transition'' in variational behavioir.   To begin the characterization 
they  need to define  wha it means to be  a ``free-time Morse minimizer''. 
 In our definition above of  `` free-time minimizer'', the minimizer $\sigma$  joined two fixed points $A, B$, and had domain a closed bounded interval $[a,b]$. 
 (Barutello, Terracini and Verzini call this type of minimizer a  ``Bolza minimizers''.)  If  the domain of $\sigma$ is the entire line $\R$
 and if    its restriction to any   compact sub-interval is a free-time minimizer in our sense, then   $\sigma$ is called  a free-time Morse minimizer. 
 (In \cite{Sanchez} this type of minimizer is called a   
 ```global free time minimizer''.)    {\bf The value $\alpha_*$ is the unique value of the homogeneity for which
 a free-time Morse minimizer exists.}  To see what happens, choose endpoints $A_N$ on the   ray $0 \xi_-$ and $B_N$ on the ray $0 \xi_+$
 with $|A_N| = |B_N| = N$ and join $A_N$ and $B_N$  by a free-time minimizer   $\gamma_N$. 
 Now let $N \to \infty$.  For $\alpha > \alpha_*$ the  curves $\gamma_N$ disappears  off to infinity.
 For $\alpha < \alpha_*$ the curves  $\gamma_N$ converge to the collision-ejection
 concatenation of two parabolic minimizers, which is the case where the Marchal lemma fails.  Exactly
 at $\alpha_*$ we get a nice convergence of $\gamma_N$ to a free-time Morse minimizer.
 This limit curve  is a parabolic non-colliision solution to Newton's equation
 connecting $\xi_-$ to $\xi_+$.  
   
To compare our set-up in Theorem \ref{JMMarchal} with theirs, observe that we have specialized to the case    $\xi_- = N$,
 $\xi_+ = S$ but    have relaxed the condition that these two points are the only two absolute minima,
thus allowing    for the minimum level set to be a continuum, as it will  be for any potential with a continuous symmetry group.
  
  \section{Back to Why}
  
 Armed with the knowledge of how the Marchal lemma can   fail, let us ammend and  return to our original   question: 
  ``Why does Marchal's lemma hold for the power law potentials?''
I will not give a full metric-inspired proof of the lemma, but rather   a sketch of plausible geometric mechanisms  behind the lemma.  
  
  The heart of the idea  is contained in section \ref{KeplerCone}.  
  I would like to   answer  ``Marchal's  lemma holds due to the conical nature of the JM metric
 near collisions.''. This answer is  incomplete
 for two reasons.  First, at total collision the metric is not actually conical, but rather is  conformal
 to a conical metric (Proposition \ref{prop_normalform}), and we have seen that certain conformal factors
 can  make the lemma fail  (Theorem \ref{JMMarchal}).   The second reason is that  there is no longer
 a single total collision point, but rather an entire collision locus and the nature of the metric
 depends on which stratum of this locus we are approaching.

Before  describing this plausibility mechanism, suppose that  I could show that Theorem \ref{JMMarchImpliesMarch} held
for the power law potentials, which is to say
that for these potentials the zero energy JM Marchal lemma implied the standard Marchal lemma.  Then it
would be legitimate to   focus my attention on  establishing
the zero energy JM Marchal lemma. 

 The zero energy JM-metric is a Riemannian metric  defined {\it away} from collisions and
so yields a metric distance function on $(\R^d)^N \setminus \Sigma$ where $\Sigma$ denotes the collision locus.  When $0 < \alpha < 2$
and $d > 1$ any point on the collision locus can be reached by any  point on the non-collision locus by a path of finite JM length.
Moreover all paths to infinity have infinite JM length. 
Thus, the metric completion of   $(\R^d)^N \setminus \Sigma$ with this distance function is all of  $(\R^d)^N$. The assertion of Marchal's lemma
now becomes {\it any geodesic for this metric which ends in collision is inextendible}.  See Appendix A for the notion of geodesics on metric spaces,
and    definition \ref{inextendible}  there for that of  a  geodesics being inextendible. 

To this end, let $\gamma$ be any minimizing geodesic ending in collision.  Then any subsegment of $\gamma$ is also a minimizing geodesic.
Moreover $\gamma$ cannot lie on the collision locus for if it did its length would be infinite.  So there is a non-collision point $\gamma(a)$ along $\gamma$.
Let $b$ be the first collision time greater than $t = a$ along $\gamma$, so that $\gamma(b) \in \Sigma$ while $\gamma([a, b))$ is collision-free.
The aim is to show that $\gamma([a,b])$ is inextendible.

If $\gamma(b)$ is a total collision point then we can use the normal form (Prop. \ref{prop_normalform})   for the metric.  Let $C$ denote the central configurations,
which is to say, the critical points of $\hat U$ on the unit sphere $S^{M-1} \subset \R^M$, $M = dN$.  For each $c \in C$,
the ray $r(t) = tc$, reparameterized, and traversed backwards, is a geodesic ending in total collision.
These rays are  the usual homothetic parabolic  central configuration solutions in the N-body problem.  Then
$\gamma$ is asymptotic to one of these  central configuration geodesics as $t \to b^-$.   Any geodesic extension of $\gamma([a,b])$
past $b$ would be asympotic to some other central configuration geodesic as $t \to b^+$.    Using the metric dilations of the cone, we
can expand to the limiting case where both the incoming and outgoing geodesics are central configuration rays.   

We are now in a situation similar to that investigated in the proof of Theorem \ref{JMMarchal}.  
There we proved the geodesic could be extended.  We want to prove our initial geodesic cannot be extended. 
Any extension would correspond to concatenating the incoming geodesic
with an outgoing one, itself   asymptotic to another  outgoing central configuration ray. 
So we reduce, as in the usual proof of Marchal, to the case of two central configuration rays concatentated at total collision.
The idea for showing that this concatenation cannot be a minimizer is to first establish  that between
any two central configurations there is a  low mountain pass. In other words, given any two
 central configurations,  there is a path on the sphere which  joins them and  along
which the maximum of $\hat U$ is fairly small, while the path's  length is   not so long.
We will call such a path a ``mountain pass path''. 
In the three-body case,  imagine the incoming and outgoing rays to be the positively and negatively
oriented Lagrange triangles -- the North and South pole of the shape sphere --  while the mountain pass curve cuts through the collinear equator at an Euler point.
The alleged result then, would be  that if we work on the cone over this mountain pass path, we can connect the
two central configurations without ever touching total collision, or indeed  any collision,
since this mountain pass path will be collision-free.

What do we do if $\gamma(b)$ is not a total collision?  
My idea here is more  vague.  I  propose copying the transformations around cluster expansions as  used
in the inductive arguments of one or the other of  the standard proofs of Marchal's lemma (\cite{Ferrario}, \cite{Venturelli}) 
and arguing that near collision the metric approximately ``splits''
into a metric  governing motion normal to the stratum (``cluster type'') of the
collision point, and a metric  tangential to the stratum, the latter being times a blowing up factor $r^{-\alpha}$ where
$r$ is distance from the stratum.  Project geodesics  onto
the normal part of the metric.  Hope that the structure of the metric along this normal part is  such that   the previous
(and incomplete)  total collision argument can be used.

This sketch of a proof is  conjectural, but the picture it gives affords me some satisfaction. 
Morevoer, the argument can legitemately  be run backwards, since we know that
Marchal implies the zero energy JM-Marchal lemma (Theorem \ref{MarchImpliesJMMarch}), thus   giving some information regarding
the geometry around mountain passes associated to the power law   $\hat U$ on the sphere.

  \eject
  
  \appendix
  
 \section{Appendix. Metric Geometry and Collision Mechanics.}
     
  We recall a few concepts from metric geometry and relate them to the JM metric. 
  A  ``length space'' is a metric space $(M, d)$ such that the distance $d(p,q)$  between any two points
 $p,q \in M$ is equal to infimum of the lengths of the paths joining the two points.  A minimizing geodesic 
 between p and q is a   curve  joining them which realizes this infimum: its   length equals $d(p,q)$. 
  Any minimizing geodesic $c$ can be parameterized by arclength, 
 in which case $d(c(s), c(t)) = |s-t|$ for all $s, t$ in the domain of $c$.  
 
 What then, is a geodesic in $M$?  Let $I \subset \R$
 be a sub-interval, possibly infinite in one or both directions.  A `geodesic'' in $M$ is a curve $c: I \to M$  which is
 parameterized by arclength and such that about any point $t_0$ in the interior of the time interval $I$
 there is an $\epsilon > 0$ such that the restriction of $c$ to $[t_0 - \epsilon, t_0 + \epsilon]$ is a minimizing geodesic between
 its endpoints. (If $t_0$ is an  endpoint of $I$ we make a similar  minimality requirement, using intervals having  one endpoint $t_0$ instead.)
 \begin{definition}
 \label{inextendible} A geodesic $c:[a, b] \to M$ is inextendible beyond  $q = c(b)$
 if it admits no extension $c:[a, b+ \epsilon) \to M$ which is a geodesic.
 \end{definition} 
 
 \begin{example} If $M$ is the upper half plane with its usual Euclidean metric,
 then the geodesics are line segments lying in $M$.  A geodesic is inextendible if
 it begins in the interior of the upper half plane and ends on the boundary.
 \end{example}

 The JM length functional defines a metric geometry away from the Hill boundary and collisions,
 which is to say, on the domain $\{ 0 < h + U < + \infty \}$ of $\R^M$.  
 The  distance $d(p,q)$ between two points $p, q$ of the domain is    the infimum of the lengths $\ell(c)$ of all paths $c$ lying in the
 domain  and  joining $p$ to $q$.    
 
 However, this  distance will typically  not complete be complete.  
 We can try to complete it by   adding  points of the closure of the domain, which is to say the Hill boundary $\{h + U =0 \}$ and  the collision locus
 $\{U = + \infty\}$, together, possibly, with some   points at infinity.    Now if   $p$ and $q$
 lie in the same path component of the Hill boundary we can travel from p to q  by a path $c$ lying in the   boundary
 in which case $\ell(c) =  d(p,q) = 0$.  So we have a choice: accept that the extended metric
 is not a true metric, but rather only a semi-metric, or collapse each boundary component of the Hill boundary to a single point.
 
 Suppose, for simplicity, that our potential is negative so that the Hill boundary is empty when the
 energy is zero or negative.   Then we can ignore problems with the Hill boundary. Suppose also that the 
 collision  locus   $\Sigma = \{ q: U(q) = + \infty \}$ does not separate $\R^M$.  Then the  distance $d(p,q)$ is well-defined and finite    for any two non-collision points.
 Let us also suppose that any path tending to infinity (meaning along which  $\|c (t) \|$ is unbounded) has infinite length.
 Then the metric completion of $\R^M \setminus \Sigma$ is obtained by adding back in all those
 collision points   for which there is a finite length path which ends in them.

For example, with the power law potentials, this completion strategy works beautifully when $0 < \alpha < 2$.  
When $\alpha = 2$ and the energy is zero, then all paths to collision have infiniite length and so
$\R^M \setminus \Sigma$ is already complete. (See \cite{pants}.)   When $\alpha > 2$ any path
ending in collision continues to have infinite length, but now  there are    finite length radial  paths
extending  out to infiniity , so that points (only one?) at infinity must be added to obtain the metric completion of $\R^M \setminus \Sigma$.

The  general case of negative  potentials having  homogeneity $0 < \alpha < 2$ which are  not  power law potentials,
and which have  collision points besides the origin, could  yield JM metrics with quite complicated   metric completion, 
depending crucially on how the potential blows up along
``spherical directions'' as  points of  the collision locus are  approached. 
It would be of interest to encounter important or useful examples of such
potentials.  

 \section{Appendix. Metric Geometry of Hill boundary}
  The Jacobi-Maupertuis  metric degenerates to zero  at the Hill boundary  $\{h_0 + U = 0\}$.
    Newtonian solutions with energy $h_0$ reflect off the  Hill boundary, retracing their path.
  At the instant $t_0$ of reflection the   kinetic energy is zero
  and hence the  velocity $\dot q (t_0)$  is zero.   
  Such solutions are called ``brake solutions''  with $t_0$ the brake instant. 

 To establish  that this retracing occurs observe that whenever $q(t)$ is a solution to Newton's equations,
  so is $q(t_0 -t)$.   Now both $q(t_0 + t)$ and $q(t_0-t)$ have that same initial conditions,
  namely $(q(t_0), 0) \in T \R^M$, at $t =0$.  It follows
  from uniqueness that the two solutions are equal:   $q(t_0 + t) = q(t_0 -t)$, which is the assertion of  retracing.
  Geodesics cannot retrace their own paths. Brake solutions,  
  insofar as they can be considered to be JM geodesics,  are inextendible past the brake moment.

  For some applications and surprises of the local geometry near the Hill boundary see 
   \cite{Seifert}, \cite{HillBoundary}.

{\bf Acknowledgements}

I would like to thank Andrea Venturelli, Mark Levi,    Alain Albouy,  Alain Chenciner,  Vivina Barutello,  Rick Moeckel, and Hector Sanchez
for helpful conversations.  I   thankfully acknowledges support by the NSF,  grant DMS-20030177.

\section{}

\end{document}